\documentclass[12pt]{article}
\usepackage{geometry}                
\usepackage{amssymb}
\usepackage{epstopdf}
\usepackage{amsmath}
\usepackage{amsthm}

\newtheorem{theorem}{Theorem}[section]
\newtheorem{lemma}[theorem]{Lemma}

\numberwithin{equation}{section}

\newcommand\RR{{\mathbb R}}

\def\CC{{\mathbb C}}
\newcommand\ep{\varepsilon}
\newcommand{\g}{\mathfrak g}

\newcommand\Pt{\tilde{P}}
\newcommand\qt{\tilde{q}}
\newcommand\pt{\tilde{p}}

\newcommand\sech{\mbox{sech}}

\newcommand{\abs}[1]{\left| #1 \right|}

\newcommand{\Rprod}[2]{\langle #1 , #2 \rangle}
\newcommand{\Lnorm}[2]{\left| \left| #2  \right| \right|_{L^{#1}}}

\begin{document}

\title{Shape of the ground state energy density of Hill's equation with nice Gaussian potential}

\author{Jos\'e A. Ram\'{\i}rez
 \and Brian Rider}


\maketitle

\begin{abstract}  Consider  Hill's operator $Q = - \frac{d^2}{dx^2} + q(x)$ in which the potential 
$q(x)$ is an almost surely continuous and  rotation invariant Gaussian
 process on the circle $x \in [0,1)$.    Viewing the classical Riccati map as a change of measure,
 we establish functional integral formulas for the
probability density function of the ground state energy and also determine the density's shape.
\end{abstract}

\section{Introduction}

Consider Hill's operator, or  $Q = - \frac{d^2}{dx^2} + q(x)$ with periodic
boundary conditions, fixed here on $S^1 = [0,1)$, 
and where the potential 
$q$ is a stationary Gaussian process.    Our aim is to characterize  the
shape of the probability density function 
of the ground state eigenvalue
$\Lambda_0(q)$ for a class of potentials.  Random Schr\"odinger operators of type 
$Q$ arise in models of disordered  solids, as is 
explained in the comprehensive texts \cite{LGP} and \cite{PF}.   There exist large bodies
of work on the bulk properties of their spectrum, 
and also on the
spectral edge which is of special importance.  The latter results are
 overwhelmingly thermodynamic in nature:  
asymptotic laws on the ground state energy  in the limit
of infinite volume (see  \cite{Sn}, \cite{Mrk}, references therein, and also \cite{M1} for a different viewpoint).  
While dimension one with Gaussian  potential presents a simplifying caricature, it is a starting point in understanding the precise
statistics of the edge in a fixed domain.

In perhaps the simplest case, when $q$ is a White Noise,
the authors and S. Cambronero  \cite{CRR} proved that the probability density function $f_{WN}$ 
of $\Lambda_0(q)$ satisfies
\begin{eqnarray}
   f_{WN}(\lambda)  =  \left\{ \begin{array}{ll} \sqrt{ \lambda / \pi} \exp{ \Bigl[ - \frac{1}{2} \lambda^2 - \frac{1}{\sqrt{2}} \lambda^{1/2} \Bigr]} 
   (1+o(1)) &  \mbox{ as  }  \lambda \rightarrow + \infty, \\
        ({4 | \lambda| }/{3 \pi }) \exp{\Bigl[ - \frac{8}{3} |\lambda|^{3/2} - \frac{1}{2} |\lambda|^{1/2} \Bigr]} 
       (1 +o(1))  &
        \mbox{ as  }  \lambda \rightarrow - \infty. \end{array} \right.
\label{WN}
\end{eqnarray}
Interestingly, the $3/2$ exponent in the left tail is shared by the analogous tail in the Tracy-Widom laws of random matrix theory.
The reason for this has recently been made clear in \cite{RRV}, which shows the Tracy-Widom distributions coincide with those
of the ground state eigenvalue for  $-d^2/dx^2 + x +$``White Noise'' on the half-line.   
In any event, and keeping with the matter at hand, it is natural to ask whether the behavior (\ref{WN}) is exhibited 
across  a rich class of random potentials.

While this exact question remains unanswered, we do provide some  progress on the general front 
by establishing the same type of result  for a large class of Gaussian potentials with continuous paths.
In particular, 
let $q$ be a stationary Gaussian process of periodicity one such that
\begin{equation}
\label{ass1}
    P \Bigl( q(\cdot) \in C[S^1] \Bigr) = 1, \ \      E[ q(x) ] = 0, \ \      E [ q(x) q(y) ]  =  K(x-y),
\end{equation} 
where the (continuous, periodic) covariance kernel $K$ satisfies 
\begin{equation}
\label{ass2}
                   \int_0^1 K(x) dx  > 0.
\end{equation}
For later use also bring in the corresponding covariance operator ${\bf
  K}f(x) = \int_0^1 K(x-y) f(y) dy$ in which $K$ 
is understood to be extended periodically.  Then (\ref{ass2}) reads ${\bf K} 1 > 0$, and we have the following.

\begin{theorem}
\label{thm:main}
For all  rotation invariant Gaussian potentials $q$ satisfying (\ref{ass1})
and (\ref{ass2}), the probability 
density function $f(\lambda)$ 
of the ground state eigenvalue $\Lambda_0(q)$ has the shape:
\begin{equation}
\label{main}
\lim_{\lambda \rightarrow +\infty} \frac{1}{\lambda^2} \log{f(\lambda)}= -\frac{1}{2 \int_0^1 K(x) dx},
\   \   \
\lim_{\lambda \rightarrow -\infty} \frac{1}{\lambda^2} \log{f(\lambda)}= -\frac{1}{2 K(0)}.
 \end{equation}

 \end{theorem}

\noindent
{\em Remark 1.} Of course, $K(x) \le K(0)$, and so the left tail is heavier than the right, as in the White Noise 
case.  This has a clear explanation: large negative deviations require only  a single excursion of the potential
while large positive deviations imply a shift of the spectral bulk.

\bigskip

\noindent
{\em Remark 2.}  For a concrete example, take the periodic Ornstein-Uhlenbeck process of mass $m > 0$, in which case  
${\bf K}^{-1} =  -d^2/dx^2  + m^2$ on $S^1$ and $K(x) = \frac{1}{2m } ( \frac{e^{mx}}{e^m -1} -\frac{e^{-mx}}{e^{-m} -1})$. Then,
the right tail decays with rate $m^2/2$ whereas the left tail has a rate of $m(\coth m - \mbox{csch} \, m)$.

\bigskip

The above should be compared with an old result of Pastur \cite{P1}.  Take $Q_L = - \triangle + q(x)$  on the cube 
of side-length $L$ in $\RR^d$  where $q$ is stationary Gaussian with covariance $K$ decaying at infinity and satisfying 
an estimate of the type $|K(x)-K(0)| \le (\log|x|)^{-\alpha}$ for $\alpha > 1$ and $x \downarrow 0$.\footnote{This will insure continuous paths, 
see \cite{AT} for example.}
Then, $\lim_{\lambda \rightarrow -\infty}  \lambda^{-2} \log N(\lambda) = -1/(2K(0))$ where $N(\lambda) = \lim_{L \uparrow \infty} (1/L^{d}) \times 
\#\{ \mbox{Dirichlet eigenvalues of } Q_L \le \lambda \}$.   Moreover, it has been pointed out to
us that the ideas behind this appraisal (which are elaborated upon in  Section 9 of \cite{PF}) allow for a simple proof of the
tail estimate of the {\em distribution} function of $\Lambda_0(-\triangle  + q)$ for continuous $q$ in finite volume.   Still, the point
here is to deal with the density directly.  Further,  our method provides the possibility of establishing
higher order asymptotics, akin to (\ref{WN}), and also provides certain insights as to the ``optimal"
 potential $q$ accounting for the tail shape, as will be explained.

Next, while
the Gaussian asymptotics in both tails is perhaps not as rich as the White Noise result,
Theorem \ref{thm:main}  raises two interesting questions. 
First,   at what regularity in the potential (or in its correlation structure)  does one see a phase transition 
between the above behavior in the left tail and the $3/2$ exponent seen for $q = $ White Noise?
Second, what happens if $ {\bf K} 1 = 0$?  While this assumption might appear to be of a technical nature, it is essential to the present method
  and, at least intuitively, the outcome.
Consider the spectral representation of $q$. Being compact and non negative definite, 
$\bf{K}$ has discrete non negative eigenvalues $\lambda_n =  \widehat{K}(n) $ associated to 
each of the trigonometric eigenfunctions ${\mathfrak e}_n(x) = e^{2\pi i n x}$, and
$
q(x) = \sum_{n= - \infty}^{\infty}  \lambda_n \g_n {\mathfrak e}_n(x), 
$
for  $\g_{n} = \bar{\g}_{-n}$  independent standard Gaussians. 
When $ {\bf K } 1 > 0$ this expansion contains the constant mode.  
The opening move of our proof hinges on the presence of this mode. 
Further, it is this mode which allows 
the potential to sit in a vicinity of $\lambda \uparrow \infty$, accounting for the Gaussian shape in that direction.

To compare (\ref{WN}) and (\ref{main}) more fully requires at least a brief sketch of the underlying methods.
The White Noise result was based on earlier work of Cambronero-McKean \cite{CM} who proved 
a functional integral expression for the density in question:
\begin{equation}
\label{CM0}
  f_{WN} (  \lambda )  =   \frac{1}{\sqrt{2 \pi}} 
   \int_{ H}   e^{- \frac{1}{2}  \int_0^1 (\lambda + p^2)^2 }    \int_0^1 e^{ 2 \int_0^x p}  \times \int_0^1 e^{-2 \int_0^x p}  \,   d P(p),
\end{equation}
where  $H$ is the space of periodic, mean-zero paths on $[0,1]$, and $P$ is a probability measure on $H$.   
More precisely, it is 
the mean-zero Circular Brownian Motion $-$ the standard Brownian Motion conditioned so that 
the Lebesgue-distributed $p(0)$ equals $p(1)$, and then further conditioned so that $\int_0^1 p = 0$.\footnote{
$P$ is also commonly referred to as the Gaussian Free Field on $S^1$.}
When $\lambda \rightarrow +\infty$,  it is clear that it is most advantageous for the path $p$ 
to concentrate in a vicinity of the origin, explaining why $f_{WN}(\lambda) \sim e^{ - \frac{1}{2} \lambda^2}$.  
When $\lambda \rightarrow - \infty$,   the factor $\int_0^1 (|\lambda| - p^2)^2$ makes the path prefer
either $+ \sqrt{|\lambda|}$ or  $- \sqrt{|\lambda|}$.  But, since $\int_0^1 p = 0$, it must split its time among both levels, and all of this
must be balanced by the Circular Brownian energy $\int_0^1 {p^{\prime}}^2$.  The outcome of this competition is the heavier 
tail.


The origin of (\ref{CM0}) is the classical Riccati transformation:  If $\psi$ solves the Hill equation $ Q \psi = \lambda \psi $, then, 
when sensible,  $p = {\psi^{\prime}}/{\psi}$ is a solution of the Riccati equation $ q = \lambda + p^{\prime} + p^2$.   
Considering the latter as a change of measure from the potential space to the new space of $p$-paths, the formula (\ref{CM0}) ensues.
For White Noise $q$,  it is clear that $p$ is locally Brownian.  Here we take on 
far smoother $q$, and so, if an analogous $p$-integral exists,  the energy term will be more coercive than
 $\int_0^1 {p^{\prime}}^2$.  This begins to explain the thinner left tail: while  $p$ may want to live
 near $\pm \sqrt{|\lambda|}$, the energy in the measure rules out the required excursion.

Actually though, Theorem \ref{thm:main} does not rely on exactly this style of integral formulation.
For the potentials considered here,  the approach behind (\ref{CM0}) runs up against certain complications.   
We use instead the following, which is also  based on the Riccati map.

\begin{theorem}
\label{thm:newdens}
Let $P$ be the law of $q$ and $\Pt$ the law of $\qt = q -\int_0^1 q$
for  $q$ as above.   With
 ${\bf K} 1 = 1$ henceforth assumed for convenience, 
\begin{equation}
f(\lambda) = \frac{1}{\sqrt{2 \pi}} \int_H e^{- \frac{1}{2} (\lambda + \Phi(\qt))^2}  d\Pt(\qt),
\label{eq:DensityFormula}
\end{equation}
where  $\Phi$ is  a positive, non-linear functional of the path, defined below in (\ref{Phi}).   Furthermore $f$ is  $C^{\infty}$ and strictly positive.
\end{theorem}

Being so 
implicit,  this formula is neither as attractive or powerful of those of Cambronero-McKean. Nevertheless, it does represent 
the ground state eigenvalue density as an expectation, and that is all that is required for the task at hand.  As already indicated,
our analysis will
show that as $\lambda \rightarrow + \infty$,  the integral  (\ref{eq:DensityFormula})   concentrates $q \sim \lambda$.  This is an in the
 White Noise case,  $p = 0$ corresponding to 
$q = 0$.  Moreover, the present result for $\lambda \rightarrow -\infty$ stems  from $q  \sim \lambda K(\cdot)/K(0)$.  That is, $q$ concentrates 
at a scaled version of the potential's covariance kernel.  Back in the White Noise case, \cite{CRR} shows
that the $3/2$ exponent for $\lambda \rightarrow -\infty$  is tied to the optimal path in (\ref{CM0}) 
being $p \sim \sqrt{|\lambda|}  {\mbox{sn}} (\sqrt{|\lambda|}  \cdot)$. Via
Riccati, this implies the optimal White Noise potential satisfies $q \sim  2 \lambda \sech^2( \sqrt{|\lambda|} \cdot)$,
an approximate Dirac measure.  In this manner, both results may be seen in a unified light.

The next section provides further background on the Cambronero-McKean approach, the difficulties it runs into in the present
set-up, and our alternate density formula (\ref{eq:DensityFormula}).  Section 3 proves  the
asymptotics for $\lambda \rightarrow +\infty$ and the lower bound as $\lambda \rightarrow -\infty$,  all of this being rather straightforward.
Similar to the White Noise calculations in \cite{CRR}, the chief difficulty  lies in the upper bound in the limit of large negative $\lambda$, and this
occupies Section 4.

\section{Integral formulas for the density}

Both are based on the Ricatti transformation. Classically, given a non-vanishing solution $\psi$ of
the Schr\"odinger equation
$  Q \psi (x) = - \psi^{\prime \prime}(x) + q(x) \psi(x) = \lambda \psi(x)$, the substitution $ p = \frac{\psi^{\prime}}{\psi}$  
produces a solution of the Ricatti equation $q(x) = p^{\prime}(x) + p^2(x) + \lambda$.   The idea in \cite{CM} was to view this 
map $q \mapsto p$ as a change in measure.


The procedure is easier to describe first for  the distribution function of $\Lambda_0$.
Fix a $\lambda$ to the left of $\Lambda_0$. Then, there is a unique positive $\psi$ solving $Q \psi = \lambda \psi$ with 
multiplier $m > 1$, {\em{i.e.}}, $\psi(x+1) = m \psi(x)$, and $ p  = {\psi^{\prime}}/{\psi}$ is a {\em periodic} solution of Ricatti's
equation with $ \log m = \int_0^1 p(x) dx > 0$.  That is, the $q$-event  $\{ \Lambda_0(q) > \lambda \}$ maps to 
the $p$-event $\{ \int_0^1 p  > 0\}$, and, with $q$ distributed according
to the Gaussian measure with covariance  operator $\bf{K}$, you will understand why
it might hold  that
\begin{eqnarray}
\label{CM1}
P( \Lambda_0(q) > \lambda ) 
& = &     \frac{1}{Z}  \int_{\{ \Lambda_0(q) > \lambda \}}     e^{-\frac{1}{2} \int_0^1 q  {\bf K}^{-1} q } \, d^{\infty} q   \\
& = & \int_{\{  \int_0^1 p >0 , \, p(0) = p(1) \}}   e^{  - \frac{1}{2} \Gamma_{\bf{K}}(p)  }   J(p) \, \times \frac{1}{Z^{\prime}}
e^{ -\frac{1}{2} \int_0^1 p^{\prime} {\bf K}^{-1} p^{\prime} } \, 
            {d^{\infty} p}.   \nonumber
\end{eqnarray}
Here
\begin{equation}
\label{factor1}
  \Gamma_{\bf K}(p) =    \int_0^1 ( \lambda + p^2) {\bf K}^{-1}  (\lambda + p^2 )  \, dx 
                                         + 2   \int_0^1  (\lambda + p^2) {\bf K}^{-1} p^{\prime}   dx       ,
\end{equation}
$Z$ and $Z^{\prime}$ are normalizers, and $J(p)$ is a Jacobian factor.   One result of  \cite{CM} is that $J(p) =  2 \sinh ( \int_0^1 p )$.

The interpretation of this expression requires some care.  When $q$ is a White Noise for example, ${\bf{K}} \equiv 1$, the reference 
measure on the right of (\ref{CM1}), or  $ \frac{1}{Z^{\prime}} e^{- \frac{1}{2} \int_0^1 {p^{\prime}}^2}  d^{\infty} p$  subject to 
$p(0)= p(1)$, can only be Circular Brownian Motion.  The first factor in (\ref{factor1}) produces the exponent $\int_0^1 (\lambda + p^2)^2$
familiar from (\ref{CM0}). The second factor, $2 \int_0^1 ( \lambda +p^2) p^{\prime}$,  which looks troublesome,  is understood to vanish by
the periodicity of the path.  This last claim is made rigorous in the course of the proof in \cite{CM},  which is effected by 
passing to the limit in appropriate discretizations of the transformation. 

To obtain a density formula  for  $\Lambda_0(q)$,  Cambronero-McKean perform a more elaborate change 
of measure based on the Ricatti correspondence, tracking the map $( q, \log m) \mapsto (p, \lambda)$.  Indeed, the potential $q$ together with the 
multiplier $m$ uniquely determine the both $p$ and the eigenvalue $\lambda$.   The upshot  is a similar integral expression,  
but restricted to the surface $\int_0^1 p = 0$ rather than the half-space $\int_0^1 p > 0$, and with a new Jacobian $J_0(p) =
  \int_0^1 e^{2 \int_0^x p} \times  \int_0^1 e^{ -2 \int_0^x p}$:
\begin{equation}
\label{CM2}
  P (\Lambda_0(q) \in d \lambda) / d \lambda = \int_{\{  \int_0^1 p = 0 , \, p(0) = p(1) \}}     e^{  - \frac{1}{2} \Gamma_{\bf{K}}(p)  }  J_0(p)  \, \times
  \frac{1}{Z^{\prime \prime}}  e^{ -\frac{1}{2} \int_0^1 p^{\prime} {\bf K}^{-1} p^{\prime} } \, 
            {d^{\infty} p}.
\end{equation}
A direct comparison with (\ref{CM0}) is now available.  Note that the presence of the conditioning, $\int_0^1 p = 0$,
requires a new normalizer denoted by $ Z^{\prime \prime}$.

Along with White Noise potential,  \cite{CM} details what happens for $q$  a periodic Ornstein-Uhlenbeck process, which, as 
we have noted before, falls into the class of potentials of present interest.  
Still, there is a serious hitch in attempting the above program for our entire class.
 Implicit in the
interpretation of (\ref{CM1}) and (\ref{CM2}) is that we can rewrite events in $q$, under the Gaussian measure $dP(q) =  \frac{1}{Z} 
   e^{-\frac{1}{2} \int_0^1 q  {\bf K}^{-1} q } \, d^{\infty} q,  
   $ in terms of $p$-expectations, with $p$
distributed according to the new Gaussian measure $ d\hat{P}(p) =  \frac{1}{Z^{\prime}} e^{ - \frac{1}{2} p^{\prime} {\bf K}^{-1} p^{\prime} } \, 
            {d^{\infty} p}$, conditioned on $p(0) = p(1)$ and perhaps $\int_0^1 p = 0$.  Unfortunately, on account of the non-linearity of the map $q \mapsto p$, it
            is possible that the law of $q$ under $P$ is not absolutely continuous with respect to its law induced by $p$ under $\hat{P}$.
This can occur when, for instance, the spectral representation of $q$ has gaps or is comprised of only finitely many modes (since then 
 the corresponding representation of $p$ will typically have more $-$ even infinitely many more $-$ modes).

The upshot is that in these ``bad" cases some or all of the expontential term $\Gamma_{\bf K}(p)$  really is troublesome
and cannot be considered as part of a density with respect to $\hat{P}$.  It must be incorporated into the measure itself $-$ 
no level of care in interpretation will save matters.  Of course, adding this non-linear object to the reference measure ruins the Gaussian nature of
the latter and makes computation difficult.  For these reasons we prefer the route described next.

\subsection*{Proof of Theorem \ref{thm:newdens}}

The idea is to carry out the Riccati transformation on only part of the path space. 
Being a product measure, we can decompose $P$ as in 
\begin{equation}
dP(q) = e^{-q_0^2/2} \frac{dq_0}{\sqrt{2\pi}} \times d\Pt(\tilde{q})
\label{eq:MeasureDec}
\end{equation}
where $q_0 =\int_0^1 q$ and $\qt = q-\int_0^1 q$ and $\Pt$ the law of $\tilde{q}$.
We similarly decompose the Riccati map.  Denote the log-multiplier by 
\[
s = \log m =  \int_0^1 p(x) dx, 
\]
where $p$ satisfies
\begin{equation}
 q = q_0+ \qt = \lambda + p' + p^2.
\label{eq:Ricatti1}
\end{equation}
Projecting onto constants, we
find that
\begin{equation}
q_0 = \lambda + \int_0^1 {\pt}^2 + s^2,
\label{eq:ProjConst}
\end{equation}
as well as
\begin{equation}
\qt = {\pt}' + {\pt}^2 - \int_0^1 {\pt}^2 +2 s {\pt},
\label{eq:ProjZeroMean}
\end{equation}
by the complementary projection.

Now fix $s$ and the potential $q = (q_0,  \qt)$.  Instead of mapping to the pair 
$(\lambda, p)$ we transform only the $q_0$ coordinate. That is, we still view  $\pt$ as a function of $q$.

\begin{lemma}
The Jacobian of the transformation $(q_0,{\qt},s)\mapsto (\lambda, {\qt},s)$  is one.
\end{lemma}

\begin{proof} 
Differentiating equation (\ref{eq:ProjConst}) in $\lambda$ one deduces that
\[
\frac{\partial q_0}{\partial \lambda}  = 1 + \int_0^1 {\pt} \frac{\partial {\pt}}{\partial \lambda},
\]
while
doing the same in equation (\ref{eq:ProjZeroMean}) shows
\[
0= \frac{\partial {\pt}^{\prime}}{\partial \lambda} + 2 {\pt}\frac{\partial {\pt}}{\partial \lambda} - 2 \int_0^1 {\pt}\frac{\partial {\pt}}{\partial \lambda} 
+ 2s \frac{\partial {\pt}}{\partial \lambda}.
\]
With $ x \mapsto \partial {\pt}/\partial \lambda (x)$ a periodic solution of this last  equation, a simple argument implies 
that in fact $\int_0^1 {\pt} \frac{\partial {\pt}}{\partial \lambda} \equiv
0$. At any zero $z$ of 
$\partial {\pt}/\partial \lambda$ we have that
\[
\frac{\partial {\pt}^{\prime}}{\partial \lambda}(z) = 2 \int_0^1 {\pt}\frac{\partial {\pt}}{\partial \lambda}, 
\]
but this is impossible unless the right hand side vanishes as it would otherwise require 
that the function $ \partial {\pt}/\partial \lambda (x)$ is, say, increasing at each of its roots.  
Hence, $\partial q_0/\partial \lambda  = 1 $.
\end{proof}

Recalling (\ref{eq:ProjConst}) we define 
\begin{equation}
\label{Phi}
 \Phi( \qt) = \int_0^1 \pt^2(x) dx 
 \end{equation}
 and write,
\begin{eqnarray*}
 P  \Bigl(   \Lambda_0(q)  \ge  \lambda \Bigr)  & = &   \frac{1}{h}
 \int_0^h    
\int_{ \{ \Lambda_0( q_0  + \qt) >  \lambda \} }  e^{- \frac{1}{2} q_0^2 } 
    \frac{d q_0}{\sqrt{2 \pi}}  \, d {\Pt} ( \qt) \, ds  \\
    & = & \frac{1}{h}   \int_0^h    \int_{ \{ \Lambda_0( \lambda^{\prime}
 + \qt) >  \lambda \} }  
e^{- \frac{1}{2} ( \lambda^{\prime} + \Phi( \qt)  + s^2)^2 } 
    \frac{d \lambda^{\prime} }{\sqrt{2 \pi}}  \, d {\Pt} ( \qt) \, ds.
\end{eqnarray*}
Now let $h \downarrow 0$, the left hand side being independent of this variable.  
At $s = 0$, we of course have $\Lambda_0 = \lambda^{\prime}$.
Hence, the event $\Lambda(\lambda^{\prime}+ \qt) \ge \lambda$ reduces to $\{ \lambda^{\prime} \ge \lambda \} \cap H$, 
and we conclude that
\begin{equation*}
 P  \Bigl(   \Lambda_0(q)  \ge  \lambda \Bigr) 
      =        \frac{1}{\sqrt{2 \pi}}  \int_{\lambda}^{\infty}   \int_H   e^{- \frac{1}{2} ( \lambda^{\prime} + \Phi( \qt) )^2 } 
                       d \lambda^{\prime}  \, d {\Pt} ( \qt).
\end{equation*}
The above expression may now be differentiated in $ \lambda$ to produce the density 
formula (\ref{eq:DensityFormula}).
(This trick of inserting an integral over the log-multiplier variable for free in step one is borrowed from \cite{CM}.)

For later use $-$ and
to be sure that the integral (\ref{eq:DensityFormula}) is sensible $-$ 
we record the following. 

\begin{lemma}
\label{lem:ContofPhi}
The functional $\qt \rightarrow \Phi(\qt)$  is bounded and continuous in  $L^\infty$. 
In particular, 
$
\abs{\Phi(\qt)} \leq  \Lnorm{\infty}{\qt}
$
and
\begin{equation}
\abs{\Phi(\qt_1)-\Phi(\qt_2)} \leq  \Lnorm{\infty}{\qt_1- \qt_2}.
\label{ineq:ContOfPhi}
\end{equation}
\end{lemma}

With this established, the stated positivity of $f(\lambda)$ is obvious.   The smoothness lies only slightly deeper.
As $ a^n e^{-\frac{1}{2} a^2}$ is uniformly bounded 
for any $n$, dominated convergence shows 
\[
    f^{\prime}(\lambda) =  - \frac{1}{\sqrt{2 \pi}}  \int_H (\lambda +
    \Phi(\qt))  
                e^{-\frac{1}{2} ( \lambda + \Phi( \qt))^2 }  d \Pt,
\]
and  $| f^{\prime} |_{L^{\infty}} < \infty$.  The argument is repeated for higher derivatives.

\begin{proof}[Proof of Lemma  \ref{lem:ContofPhi}]
Note that the first inequality follows from the second one since $\Phi(0) = 0$.  Next, when
$s=0$, we have $q_0 = \lambda + \int_0^1 {\pt}^2 $ and ${\qt} = {\pt}^{\prime} + 
{\pt}^2 - \int_0^1 {\pt}^2 $. Let  $p_1$ and $p_2$ be  solutions of the Ricatti equation 
corresponding to  $q_1$ and $q_2$ and a common $\lambda$.  With $\Delta p = p_1 - p_2$,
{\em etc}, it holds 
\[
    (\Delta {\pt})^{\prime}(x) =  \Delta \qt (x) -  (\pt_1(x) + \pt_2(x) )  \Delta \pt (x)  +   \int_0^1 ( {\pt}_1 - {\pt}_2 )  .
\]
And so,
\[
\Delta \pt (x) = \Delta \pt(0) e^{-\int_0^x(\pt_1+ \pt_2)} + \int_0^x
e^{-\int_y^x (\pt_1+ \pt_2)} \left( \int_0^1 (\pt^2_1- \pt^2_2) + \Delta \qt (y)  \right) dy. 
\]
By periodicity we obtain 
\[
\int_0^1 (\pt^2_1- \pt^2_2)  =  \Phi(\qt_1) -  \Phi(\qt_2) 
= ( \int_0^1 e^{-\int_y^1(\pt_1+ \pt_2)} )^{-1}  \int_0^1 e^{-\int_y^1( \pt_1+ \pt_2)} \Delta \qt(y)  \, dy, 
\]
from which the bound (\ref{ineq:ContOfPhi}) follows.
\end{proof}

\subsection*{A distribution formula}

For completeness, we provide a distribution function formula for the ground state energy in the 
same spirit as our density formula.   

\begin{theorem}
\label{lem:dist}
The more restrictive map $(q_0,{\qt})\mapsto (s,{\qt})$  from 
$\{ \Lambda_0(q) > \lambda \}$ to $\{ s > 0 \}$ will produce
\begin{equation}
\label{distribution}
P \Bigl(\Lambda_0(q)>\lambda \Bigr) =   \frac{1}{\sqrt{2 \pi}}  \int_H 
 \int_0^{\infty} e^{-\frac{1}{2}\left(\lambda + \Phi(\qt) + s^2 \right)^2} J(s,\qt) \, {ds} \, d \Pt(\qt),
\end{equation}
with Jacobian factor given by
\[
J(s,\qt)^{-1} = \frac{1}{1-e^{-2s}} \int_0^1 \int_0^1 e^{2 \int_0^x \pt + s}
 e^{ - 2 \int_0^{y}  \pt + s} -  \int_0^1  \int_{x}^1 e^{2 \int_0^x \pt + s}
 e^{ - 2 \int_0^{y}  \pt + s}.
\]
\end{theorem}

Note that  $J(s, \qt)$ is positive whenever $s >0$. 
Similar to the Cambronero-McKean formulas we have here an integral
over $s > 0$, while the density integral lives on  $s=0$.

\begin{proof}
We need only the Jacobian ${\partial q_0}/{\partial s}$. 
Differentiating (\ref{eq:ProjConst}) in $s$ (rather than $\lambda$) gives
\[
\frac{\partial q_0}{\partial s} =  2 \int_0^1 \frac{\partial {\pt}}{\partial s} {\pt} + 2s. 
\]
Making the definition $A = \int_0^1 \frac{\partial {\pt}}{\partial s} {\pt}$, and as before now differentiating 
(\ref{eq:ProjZeroMean}) in the same variable produces 
\[
0 = \frac{\partial {\pt}'}{\partial s} + 2 {\pt}\frac{\partial {\pt}}{\partial s} -2 A+2 s\frac{\partial {\pt}}{\partial s}+2 \pt.
\]
This may be solved:
\begin{equation}
\label{Jeq}
\frac{\partial {\pt}}{\partial s}(x) = \frac{\partial {\pt}}{\partial s}(0) 
e^{-2\int_0^x( \pt+s)} + 2\int_0^x (A - \pt(y)) e^{-2\int_y^x (\pt+s)} \, dy,
\end{equation}
with the initial condition obtained from the periodicity. In particular,
\[
\frac{\partial \tilde{p}}{\partial s}(0) =\frac{\partial
  \tilde{p}}{\partial s}(1)  
= \frac{\partial \tilde{p}}{\partial s}(0) e^{-2\int_0^1(\pt+s)} + 2\int_0^1 (A- \pt(x)) e^{-2\int_x^1(\pt+s)} \, dx,
\]
and so
\begin{eqnarray}
\label{Init}
\frac{\partial {\pt}}{\partial s}(0)  
& =  & \frac{2 (A+s)}{ 1-e^{-2s}}  \, \int_0^1 e^{-2\int_x^1(\pt+s)} \, dx   -1 .  \nonumber
\end{eqnarray}
Next, multiplying equation (\ref{Jeq}) by $\pt(x)$ and integrating from $0 $ to $1$ provides an equation for $A$,
\begin{eqnarray*}
A  & =  &  \frac{\partial \pt}{\partial s}(0)   \left(      \frac{1 - e^{-2s}}{2} - s \int_0^1 e^{-2 \int_0^x (\pt + s)}  dx \right) 
     + 2 \int_0^1 \pt(x) \int_0^x ( A - \pt(y) ) e^{-2 \int_y^x (\pt + s)}  dy dx , 
\end{eqnarray*}
which, after substituting in (\ref{Init}) and some manipulations, simplifies to
\[ 
  1 =  \frac{2 (A+s) }{1 - e^{-2s}}  \int_0^1 \int_0^1 e^{- 2 \int_0^x (\pt + s) }  e^{ 2 \int_y^1 (\pt + s) } 
             + 2 (A+s) \int_0^1 \int_y^1 e^{- 2 \int_y^x (\pt + s)}.
\]
From here the expression for $J =  {\partial q_0}/{\partial s} = 2 A + 2s$ is easily read off.
\end{proof}

\section{Three out of the four bounds}

The asymptotics for $\lambda \rightarrow + \infty$ rest on little more than the non-negativity of
$\Phi$ and the positivity of $f$.  By the first fact, we have at once that
\[
   f(\lambda) = \frac{1}{\sqrt{2 \pi}} 
     \int_H e^{-\frac{1}{2} ( \lambda + \Phi( \qt))^2 }  d \Pt (\qt)  \le  \frac{1}{\sqrt{2 \pi}} e^{-\frac{1}{2} \lambda^2},
\]
for all $\lambda \ge 0$.  For the lower bound in the same direction, with $\tilde{E}$ denoting $\Pt$-expectations,
\begin{eqnarray*}
 f(\lambda) & =     &  \frac{1}{\sqrt{2 \pi}}   e^{-\frac{1}{2} \lambda^2} 
          \int_H e^{ - \lambda \Phi(\qt)  - \frac{1}{2} \Phi^2( \qt) }  d \Pt (\qt) \\ 
           & \ge &   e^{-\frac{1}{2} \lambda^2}    
           f(0)  \exp \Bigl[  - \lambda \tilde{E} [   \Phi(\qt)
          e^{-\frac{1}{2} \Phi^2(\qt)} ]  /   
                      \sqrt{2 \pi} f(0)    \Bigr] \\
         & \ge &       f(0)  e^{-\frac{1}{2} \lambda^2  -   \frac{1}{  \sqrt{2 \pi e} f(0)}  \lambda},                               
\end{eqnarray*} 
since $ a e^{-\frac{1}{2} a^2} \le \sqrt{e}$ for $a \ge 0$.

\subsection*{Lower bound for $\lambda \rightarrow - \infty$}

Considerations in the next section will indicate why the $f$-integral concentrates in
a vicinity of (the mean zero version of) $q_{\lambda}(x) = \lambda K(x) / K(0)$ as $\lambda \rightarrow -\infty$.   For the 
lower bound we then proceed by restricting the expectation as in 
\[
   f(\lambda) \ge  \frac{1}{\sqrt{2 \pi}}  \tilde{E} \Bigl[  e^{-\frac{1}{2} ( \lambda + \Phi( \tilde{q}) )^2} , \, ||  \tilde{q} - \tilde{q_{\lambda}} {||}_{L^{\infty}} < \ep |\lambda| \Bigr]
\]
with whatever $\ep > 0$ and $ \tilde{q}_{\lambda}(x) = \lambda  \Bigl(\frac{ K(x) - 1}{K(0)} \Bigr)  := \lambda  q_{K}(x)$.
By Lemma  \ref{lem:ContofPhi} we may then write 
\begin{eqnarray}
\label{largelower}
   \log  f(\lambda)   & \ge &   \log \tilde{E} \Bigl[  e^{-\frac{1}{2} ( |\lambda + \Phi( \tilde{q_{\lambda}})| + \epsilon|\lambda| )^2} , \, ||  \tilde{q} - \tilde{q_{\lambda}} {||}_{L^{\infty}} < \ep |\lambda| \Bigr] - \frac{1}{2}\log (2 \pi) \nonumber \\ 
&\ge&   - (\frac{1}{2} + \ep) \Bigl(  \lambda + \Phi(\tilde{q}_{\lambda}) \Bigr)^2   - (\ep+\frac{1}{2}\ep^2 ) \lambda^2  
                                 \\ 
                                  &&   \  \  \ +  \log \Pt^{(\lambda)} \Bigl(  ||  \tilde{q} -   q_K    {||}_{L^{\infty}} < \ep   \Bigr)  - \frac{1}{2}\log (2 \pi), 
                                  \nonumber
\end{eqnarray}
where $ \Pt^{(\lambda)}$ indicates the scaled Gaussian measure induced by $\Pt$ on ${|\lambda|^{-1}} \qt $.  

For the third term we immediately have that
\begin{equation}
\label{largelower2}
   \liminf_{|\lambda| \rightarrow \infty}  \frac{1}{\lambda^2} \log \Pt^{(\lambda)} \Bigl(  ||  \tilde{q} -   q_K    {||}_{L^{\infty}} <  \ep   \Bigr) 
     \ge  - \frac{1}{2}   \langle q_K,  {\bf K}^{-1} q_K  \rangle =  \frac{1-K(0)}{2K^2(0)},
\end{equation}
in which we note that with ${\bf K} 1 = 1 $ assumed, $K(0) \ge 1$ and the right hand side is $\le 0$.  This is the 
standard Schilder's Theorem 
for Gaussian processes (\cite{DS} may be consulted), and it is here that we
require the continuity of the potential sample paths. 

As for the first term of (\ref{largelower}), a bit more information on the
functional $\Phi$ is needed.   Returning to the defining relation 
(\ref{eq:ProjConst}), note that if $q$ is a Hill potential with ground
state eigenvalue $\lambda$ , 
then $\Phi(\tilde{q}) =   q_0 - \lambda$. 
That $q_0 - \lambda$ is invariant under shifting $q$ by a constant  is
another way to see that $\Phi$ is 
well defined on the mean-zero $\tilde{q}$.
Now, $q_{\lambda}(x) = \lambda K(x) / K(0)  \ge \lambda $ is {\em not} a
Hill potential for ground state energy 
$\lambda$ (unless of course
$K$ is the constant function).  Denote by $\Lambda_0(\lambda)$  the actual minimal eigenvalue for the (sure) Hill operator 
$-\frac{d^2}{dx^2} +  q_{\lambda}(x)$.  Then,    $\Phi( \tilde{q}_{\lambda}) =   1/ {K(0)} - \Lambda_0(\lambda)$,
and
\begin{equation}
\label{largelower3}
    | \lambda - \Lambda_0(\lambda)|  = o(|\lambda|)  \mbox{ will imply } 
   \liminf_{\lambda \rightarrow -\infty}  \frac{1}{\lambda^2}  
    \Bigl( \lambda + \Phi(\tilde{q}_{\lambda}) \Bigr)^2 \ge  \frac{1}{K^2(0)}.
\end{equation}
Combined with (\ref{largelower2}), and taking $\ep \downarrow 0$ after the fact, brings out the 
desired conclusion  $\liminf_{\lambda \rightarrow -\infty}$  $  \lambda^{-2} \log f(\lambda) $  $\ge - 1/(2 K(0))$.

To verify the first part of (\ref{largelower3}), by the variational characterization of $\Lambda_0 (\lambda)$ 
it trivially holds
\[
    \lambda \le  \Lambda_0 (\lambda) \le \int_0^1 (\psi^{\prime}(x))^2 \, dx  + \lambda \int_0^1  \frac{K(x)}{K(0)} \psi^2(x)  \, dx
\]  
for any choice of periodic $\psi$ such that $\int_0^1 \psi^2 = 1$.  So
choose $\psi$ to be a smooth 
approximation of the identity supported
in a neighborhood of $ x= 0$ of width $\delta > 0$.    Then, $ \int_0^1 (\psi^{\prime})^2 = O(\delta^{-1})$ while the 
difference $| \lambda - \lambda \int_0^1 (K(x)/K(0)) \psi^2(x)| \le$
$|\lambda| \int_0^1 | K(0) - K(x)| \psi^2(x)$  may be bounded by a constant multiple of 
$|\lambda| \times \sup_{|x| < \delta/2} | K(x) - K(0)|$.   As $K$ is
continuous  and $\delta$ is chosen at will, 
the proof is complete.

\section{Upper bound for $\lambda \rightarrow - \infty$}

For the $\lambda \rightarrow - \infty$ case we must at last consider the full exponent of the integral
form of the density.  That is, the asymptotics will stem from minimizing
\begin{equation}
 J_{\lambda}(\qt) = 
     \frac{1}{2} ( \lambda + \Phi(\qt) )^2  + \frac{1}{2} \Rprod{ \qt }{{\bf K}^{-1} \qt }
\label{J0}
\end{equation}
over smooth, periodic, mean-zero $\qt$ subject to $\int_0^1 p = 0$ through the Riccati map.  
 First we consider this
minimization problem in the limit of large negative $\lambda$, then we turn to the
corresponding upper bound on $f(\lambda)$.

\subsection*{The variational problem}

Minimizing  $J_{\lambda}$ is equivalent to doing the same for
\begin{equation}
\label{J}
J(q) = \frac{1}{2} \Rprod{q}{{\bf K}^{-1}q}
\end{equation}
over the set 
\[
A_{\lambda}  = \Bigl\{  q : q = \lambda + p' + p^2  \mbox{ for some }  p \in C^{\infty}   \mbox{  with } \int_0^1 p =0
                                          \Bigr\} .
\]
The needed result is the following.

\begin{theorem}  
It holds,
\label{thm:Lower}
\[
\liminf_{\lambda \rightarrow - \infty} \frac{1}{\lambda^2} 
 \left\{  \inf_{ q \in A_{\lambda}}  J(q) \right\}
     \ge  \frac{1}{2 K(0)}.
\]
\end{theorem}

The proof requires two preliminary steps.

\begin{lemma}
\label{lemma:Exists}
 For fixed $\lambda$, 
the minimum of $J(q)$ over $A_{\lambda}$ is attained at a continuous
$q_{\lambda} < 0$ which
satisfies
\begin{equation}
\label{EulerEquations}
    q_{\lambda}(x) = {\bf K} a(x),    \   \     \    \  a^{\prime}(x) = 2 p(x) a(x),
\end{equation}
with a mean-zero $p$ such $q_{\lambda} = \lambda + p^{\prime}(x) + p^{2}(x)$. 
\end{lemma}

\begin{proof}  For existence, it is simpler to take the initial viewpoint
  (\ref{J0}). For any minimizing sequence $( \pt_n , \qt_n)$,  the  $
  \Rprod{\qt_n}{{\bf K}^{-1}\qt_n} < \infty$ term and Rellich's criteria
  will show that there is a $\qt_{\infty}$ (eventually equaling $q_{\lambda} -
  \int_0^1 q_{\lambda}$) with $\qt_n \rightarrow
  \qt_{\infty}$ in $L^2[S^1]$.   The first term shows that $\Phi(\qt_n) =
  \int_0^1 \pt_n^2$ is uniformly bounded, and so we at least have a
  $L^2$-subsequence  $ \pt_{n'}  \rightarrow \pt_{\infty}$ giving a candidate pair
  to
satisfy Riccati.  Along the way: assuming $\pt_n(0) = 0$, 
\begin{equation}
\label{intric}
   \pt_n(x) = \int_0^x \Bigl( \qt_n - \pt_n^2 + \int_0^1 \pt_n^2 \Bigr), 
\end{equation}
and so $\{ \pt_n \}$ is also uniformly bounded in $L^{\infty}$ and hence is 
equicontinuous.  Then taking pointwise limits (\ref{intric}) we find that 
$\pt_{\infty}$ is absolutely continuous, and  $q_{\infty}  =
\lambda + p' + p^2$ holds almost everywhere.

Next, given existence, 
any minimum must be  a stationary point of the extended functional
\[
{J}(q,p;a,\alpha) = \frac{1}{2} \Rprod{q}{{\bf K}^{-1}q} - \int_0^1 a(x) \Bigl(q(x)-(\lambda + p^{\prime}(x)  
 + p^2(x) ) \Bigr) dx  - \alpha \int_0^1 p(x) dx,  
\]
with Lagrange multipliers $\alpha \in \RR$
and (periodic) $a \in L^2$. Upon differentiating 
we obtain the set of equations:
\begin{eqnarray}
q_{\lambda}(x)  & = & \lambda + p^{\prime}(x) + p^2(x) ,  \label{eq:Ric} \\
a^{\prime}(x) & = & 2 p(x) a(x) + \alpha,   \label{eq:ExpP} \\
q_{\lambda}(x) & = &  {\bf K} a(x),  \label{eq:RelationQA}
\end{eqnarray}
subject still to $\int_0^1 p = 0$.  Both $p$ and $a$ certainly depend on $\lambda$,
but we suppress this for now.

By (\ref{eq:RelationQA}) and the continuity of $K$, $q_{\lambda}$ is
actually continuous and we recover the  Riccati correspondence in the naive sense.
Thus, any minimizing $q_{\lambda}$ is a Hill potential with groundstate eigenvalue $\lambda$.
Focussing here on the  case
$\lambda < 0$, $q_{\lambda}$ must be negative at some point.  It follows from 
(\ref{eq:RelationQA}) and ${\bf{K}} 1 = 1$ that $a$ must also be negative somewhere.
In fact $ a(x) < 0 $ for all $ x \in S^1$.  In particular, if $\alpha = 0$ then
$a$ can be solved for explicitly from (\ref{eq:ExpP}) and is seen to be of one sign.
While, 
if  $\alpha \neq 0$,  then  $a^{\prime}(z) = \alpha$  at every $z$ such that $a(z)= 0$, but this
is impossible as $a$ is periodic.  By appealing once again to
(\ref{eq:RelationQA}) we may conclude that  $q_{\lambda}$ is everywhere negative as well. 

It is now possible to divide 
(\ref{eq:ExpP}) through by $a < 0$ and integrate from 0 to 1 to  show that
 $\alpha$ must actually be zero.   (We remark that this observation was missing in  
 \cite{CRR}, where the above formulation of the variational problem 
 was not used and additional effort  was needed to get around this point.)
In short, equation (\ref{eq:ExpP}) may be replaced with
the relation $a^{\prime}(x) =2p(x) a(x)$
advertised in the statement.
\end{proof}

In many ways the next lemma summarizes 
the difference between the White Noise case  
and the present situation of 
continuous covariance potentials.

\begin{lemma}
There is the bound 
\[
\Lnorm{\infty}{f}^2 \leq C\Rprod{f}{{\bf K}^{-1}f}
\]
for some constant $C = C(K)>0$.
\label{lemma:IneqInfK}
\end{lemma}

\begin{proof}
We show the equivalent, $\Lnorm{\infty}{{\bf K}g}^2 \leq C\Rprod{g}{{\bf K}g}$.  It holds,
\begin{eqnarray*}
\Lnorm{\infty}{{\bf K}g} &=& \sup_{\phi \in L^2:\Lnorm{1}{\phi}\leq 1} \int_0^1 ({\bf K}g)(x)  \phi(x) dx \\
& \leq & \left(\int_0^1 g {\bf K}g \right)^{1/2} \sup_{\phi \in L^2:\Lnorm{1}{\phi}\leq 1} \left( \int_0^1 \phi {\bf K} \phi \right)^{1/2} \\
& = & \Rprod{g}{{\bf K}g}^{1/2} \sup_{\phi \in L^2:\Lnorm{1}{\phi}\leq 1}  \left(\int_0^1 \phi(x) \int_0^1 K(x-y) \phi(y) dy  dx \right)^{1/2} \\
& \leq & ||K {||}_{L^{\infty}}^{1/2} \, \Rprod{g}{{\bf K}g}^{1/2}.
\end{eqnarray*}
Line two uses the fact that ${\bf{K}}$ is non-negative and symmetric.  
For the last line,  $||K{||}_{L^{\infty}} $ is finite given the continuity of $K$.
\end{proof}

\begin{proof}[Proof of Theorem \ref{thm:Lower}]
The first step is to note that the normalized family of minimizers $\{  \frac{1}{|\lambda|} q_{\lambda} \}$ is bounded 
in $L^{\infty}$ and has modulus of continuity independent of $\lambda \rightarrow -\infty$.  
From the Euler-Lagrange equations (\ref{EulerEquations}) we find that:
with $\Delta_K$ denoting the modulus of continuity of the kernel $K$,
\begin{eqnarray}
\label{Kdelta}
\abs{ {q_{\lambda}}(x)- {q_{\lambda}}(x^{\prime})} &=& 
\abs{\int_0^{1} \left[ K\left( {x-y} \right)-K\left( {x^{\prime}-y} \right)\right]  {a}(y)\, dy }  \\
& \leq & \sup_{y \in S^1} \abs{K\left(  {y+x} \right) - 
K \left( {y + x^{\prime}} \right)}  \int_0^{1} 
|   {a}(y) | \, dy  \nonumber \\
& = & \sup_{y \in S^1}  \abs{K\left(  {y+x} \right)-K\left(  {y+x^{\prime}} \right)} 
  \int_0^{1} 
 | {q_{\lambda}}(y) | \, dy  \nonumber\\
& \leq & \Delta_{K}(x-x^{\prime}) \, || q_{\lambda} {||}_{L^{\infty}}.  \nonumber
\end{eqnarray}
In line three we  have used the identity
$
     \int_0^{1}{a}(x)\, dx =    \int_0^{1}  {q_{\lambda}}(x)\, dx,  
$
implied by  $q_{\lambda} = {\bf{K}} a$ ,  ${\bf{K}}1 = 1$,
and the fact both $a$ and $q_{\lambda}$ are everywhere negative, as
is pointed out in the proof of Lemma \ref{lemma:Exists}.  Next, by 
Lemma \ref{lemma:IneqInfK} we also see that
\begin{equation}
  \frac{1}{2} \lambda^2 = J(\lambda) \geq \inf_{q \in A_{\lambda}}  J(q) = J(q_{\lambda}) 
 =  \frac{1}{2} {\Rprod{q_{\lambda}}{{\bf K}^{-1} q_{\lambda}}} \geq \frac{1}{2C} \Lnorm{\infty}{q_{\lambda}}^2 
\label{norm_bound}
\end{equation}
at any minimizer $q_{\lambda}$.  Hereafter setting $\hat{q}_{\lambda} = \frac{1}{|\lambda|} q_{\lambda}$,
(\ref{Kdelta}) and (\ref{norm_bound}) may be summarized as in 
$ \sup_{\lambda < 0}  \Lnorm{\infty}{ \hat{q}_{\lambda}} \le \sqrt{C(K)}$ and 
$  \sup_{\lambda < 0}  \Delta_{ \hat{q}_{\lambda}} \le \sqrt{ C(K)} \Delta_K$.

Next, having the equicontinuity of $\{ \hat{q}_{\lambda} \}$ and the 
fact $ {\liminf}_{\lambda \rightarrow -\infty} \{ \frac{1}{\lambda^2}
\inf_{q \in A_{\lambda} }J(q) \}$  $  = {\liminf}_{\lambda \rightarrow
  -\infty}  \{ J( \hat{q}_{\lambda} ) \}$,  
we introduce the scaled equations
\begin{equation}
\label{lastequations}
    \hat{q}_{\lambda}(x) = -1 + \frac{1}{\sqrt{|\lambda|}}  \hat{p}_{\lambda}^{\prime}(x)  +  \hat{p}_{\lambda}^2(x),  \  \ 
    \hat{q}_{\lambda}  = {\bf K} \hat{a}_{\lambda},   \   \
    \hat{a}_{\lambda}^{\prime}(x) 
           = 2 \sqrt{|\lambda}| \hat{p_{\lambda}}(x)  \hat{a}_{\lambda}(x) 
\end{equation}
where $\hat{a}_{\lambda}(x) = \frac{1}{|\lambda|} a(x) =  \frac{1}{|\lambda|} a_{\lambda}(x)$  and
$ \hat{p}_{\lambda}(x) = \frac{1}{\sqrt{|\lambda|}} p(x) =  \frac{1}{\sqrt{|\lambda|}} p_{\lambda}(x)$.  By rotation invariance
and the fact $\int_0^1 \hat{p}_{\lambda} = 0$ we may assume that
\begin{equation}
\label{center}
     \hat{p}_{\lambda}(0)  = 0,   \   \   \    \hat{p}_{\lambda}^{\prime}(0)   \le 0
\end{equation}
holds in (\ref{lastequations}) along the sequence $\lambda \rightarrow - \infty$.   Combining the second two equations of 
(\ref{lastequations}) we may also write
\begin{equation}
     \hat{q}_{\lambda}(x)  = \int_0^1  \hat{q}_{\lambda}  \times \int_0^1 K(x - y) \psi_{\lambda}(y) \, dy,
\end{equation}
in which $\psi_{\lambda}(x) =  e^{2 \sqrt{|\lambda|} \int_0^x
  \hat{p}_{\lambda} }  
/  (  \int_0^1  e^{2 \sqrt{|\lambda|} \int_0^{x'} \hat{p}_{\lambda} }  \,dx'  )$.

Now consider a minimizing sequence $\{ \hat{q}_{\lambda'} \}$.  By the
pre-compactness of both $\{ \hat{q}_{\lambda} \}$ 
and the  family of probability measures $\{ \psi_{\lambda} \}$ on $S^1$, we may choose a further subsequence along which 
\[
     \hat{q}_{\lambda''}(x)  \rightarrow  \hat{q}_{\infty}(x) \mbox{ pointwise, }    \mbox{ and } 
      \psi_{\lambda''}(x) dx  \rightarrow   \psi_{\infty}(dx)  \mbox{ weakly.}
\]
By the continuity of $K$, 
passing along this subsequence it is found that  $\hat{q}_{\infty}(x) = c_{\infty}  \int_0^1 K(x-y) \psi_{\infty}(dy)$.  In particular, 
bounded convergence implies that $ \int_0^1  \hat{q}_{\lambda''}$ settles down to some 
constant $c_{\infty}$.    To pin down $c_{\infty}$, we first note that anywhere $ \hat{p}_{\lambda}^{\prime}(x)   \le 0$ it must be that
$ \frac{1}{\sqrt{|\lambda|}} \hat{p}_{\lambda}^{\prime}(x)  \rightarrow 0$.

Assume to the contrary that say $ \limsup  \frac{1}{\sqrt{|\lambda|}}
\hat{p}_{\lambda}^{\prime}(0) \le - \delta < 0$. If to the left and right
of $x = 0$ this ratio were $\ge 0$, it would follow that $\hat{q}_{\infty}$
was discontinuous at this point, and that cannot be.  So now try with $
\hat{p}_{\lambda}^{\prime}(x) /\sqrt{|\lambda|}$ staying negative on some 
neighborhood containing the origin.  This will lead to the false conclusion 
that $\hat{q}_{\infty}$ is unbounded. 
In light of this and (\ref{center}), we 
now see that $\hat{q}_{\lambda''}(0)$ can only converge to $-1$.  
Hence, 
\[
      \hat{q}_{\infty}(x) =   - \frac{ \int_0^1 K(x-y) \psi_{\infty}(dy) }{ \int_0^1 K(-y) \psi_{\infty}(dy)} .
\]
 
To finish, by the lower semi-continuity of $q \mapsto J(q)$, we have that
\[
    \liminf_{\lambda \rightarrow -\infty} J( \hat{q}_{\lambda} )  \ge  J(  \hat{q}_{\infty}) 
     =  \frac{ \int_0^1 \int_0^1 K(x-y)  \psi_{\infty}(dx) \psi_{\infty}(dy) }{ 2  \left(   \int_0^1 K(-y) \psi_{\infty}(dy) \right)^2}, 
\] 
and it remains to see that the final ratio lies above $1/(2 K(0))$ for whatever probability measure $\psi_{\infty}$.
However, for any measures $\mu$ and $\nu$ of mass one
$
 \Rprod{ \mu}{ {\bf K} \nu }^2 \le \Rprod{\mu}{ {\bf K}  \mu } \Rprod{ \nu}{ {\bf K}
   \nu},
$
and the desired inequality holds by setting $\mu = \psi_{\infty}$ and $\nu = $ the dirac delta measure at the
origin.  Afterward, it is understood that it is optimal to have had $\mu_{\infty} = \delta_0$ and so 
$\hat{q}_{\infty}(x) =   K(x) /K(0)$, or that $q_{\lambda} \sim \lambda K(x)/K(0)$.
\end{proof}

\subsection*{Proof of the final bound}

The basic idea is similar to that behind the main Large Deviation estimate
in \cite{CRR}, though a few steps  
have been streamlined with the aid of experience.  

First,  introduce  the shorthand  $R(\lambda,\qt) = \left(\lambda + \Phi(\tilde{q}) \right)^2$.
To prove that
\begin{equation}
\limsup_{\lambda \rightarrow -\infty} \frac{1}{\lambda^{2}}   \log \int_H e^{-\frac{R(\lambda,\qt)}{2}} d\Pt (\qt)  \leq -\frac{1}{2K(0)}
\label{LDresult}
\end{equation}
we define the sets 
\[
H(\gamma,\eta) = \Bigl\{ \qt \in H :  \abs{R(\mu,\qt) - \gamma \lambda^2}\leq \eta \lambda^2 \Bigr\}
\]
for $\gamma$ and $\eta$ positive, 
and decompose  the integral on the left of (\ref{LDresult}) as in
\begin{equation}
\int_H e^{-\frac{R(\lambda,\qt)}{2}} d\Pt(\qt) 
\leq \sum_{0\leq k \leq 1/K(0)\eta} \int_{H(\eta k,\eta)}
e^{-\frac{R(\lambda,\qt)}{2}}\,d\Pt(\qt) + \exp{\Bigl[-\frac{\lambda^2}{2 K(0)} \Bigr]}.
\label{ineq:Initial}
\end{equation}
Each of the terms in the sum are in turn bounded by 
\begin{equation}
\int_{H(\gamma,\eta)} e^{-\frac{R(\lambda,\qt)}{2}} d\Pt \leq \exp\left[(-\gamma+\eta)\frac{\lambda^2}{2} \right]  
 \Pt  \Bigl( H(\gamma,\eta) \Bigr),
\label{Cheby}
\end{equation}
and to continue we further define
\[
D(\gamma)= \Bigl\{  q \in H : R(\lambda,\qt) \leq \gamma \lambda^2 \Bigr\}.
\]
Notice that $\Pt  ( H(\gamma,\eta) ) \leq \Pt ( D(\gamma+\eta)) $. 

In order to estimate
the $\Pt$-probability of the event $D(\gamma + \eta)$,  truncate the path 
$\qt$ based on its Fourier expansion. If $\qt(x) = \sum_{-\infty}^{\infty} c_k e^{ikx}$  (with $c_0 = 0$),
let $\qt_{n}(x) = \sum_{-n}^n c_k e^{ikx}$ denote the projection onto the
low modes.  Then,
\begin{eqnarray}
\Pt \Big( \qt \in D(\gamma+\eta) \Bigr)  &\leq&  \Pt  \Bigl( \Lnorm{\infty}{\qt} \geq M |\lambda |  \Bigr) + 
                 \Pt \Bigl(  \Lnorm{\infty}{\qt_{n}} \geq M | \lambda|  \Bigr)    \\
&    &+ \,   \Pt \Bigl(  \abs{R(\qt_{n},\lambda)-R(\qt, \lambda)} \geq \eta \lambda^2, \Lnorm{\infty}{\qt}\leq M | \lambda | , 
                                 \Lnorm{\infty}{\qt_{n}} \leq M  | \lambda |  \Bigr) \nonumber \\
&   & +  \, \Pt  \Bigl( \qt_{n} \in D(\gamma+2\eta) \Bigr), \nonumber
\label{ineq:DecompProbD}
\end{eqnarray}
for any positive constant  $M$. We now treat each of the terms on the right of (\ref{ineq:DecompProbD}) in turn.


\bigskip
\noindent
{\em Terms 1 and 2}: For the first, we recall 
Borell's inequality  (see  \cite{AT}, Theorem 2.1.1): for any $c > 0$,
\[
\Pt  \Bigl( \Lnorm{\infty}{\qt} - m_{\infty} \geq c \Bigr) \leq \exp{ \Bigl[-  \frac{c^2}{2\sigma^2_{\infty}} \Bigr]},
\]
where $m_{\infty} = \int \sup_x \abs{\qt(x)} d\Pt$ and $\sigma^2_{\infty} = \sup_x \int \qt^2(x) d\Pt$. In the present case,
standard facts on Gaussian processes, summarized in say \cite{AT} Theorem 1.3.3, will imply 
that $  m_{\infty} <\infty$.  Furthermore, by stationarity  
$
\sigma^2_{\infty} = \sum_{k=1}^{\infty} \lambda_k^2  < \infty,
$
the kernel $K$ being continuous (the $\lambda_k$'s are the eigenvalues of ${\bf K}$) .
Hence,
\begin{equation}
\Pt \Bigl( \Lnorm{\infty}{\qt} \geq M |\lambda| \Bigr) 
\leq C \exp{\left[- \frac{M^2 \lambda^2}{2 \sigma^2_{\infty}} \right] },
\label{ineq:No1}
\end{equation}
for some constant $C$ and 
$|\lambda|$  large enough.  It is plain that $ \Pt(  \Lnorm{\infty}{\qt_n}
\geq M |\lambda |  )$ 
satisfies an identical bound.

\bigskip

\noindent
{\em Term 3:} Since for any $f$ and $g$ $\in H$,
\[
\abs{R(\lambda,f)-R(\lambda,g)} = 2 | \lambda| \abs{ \Phi(f) - \Phi(g)} + \abs{\Phi^2(f) - \Phi^2(g)},
\]
Lemma \ref{lem:ContofPhi} implies that
\begin{equation*}
\label{Rdiff}
\abs{R(\lambda,f)-R(\lambda,g)} \leq ( 2 |\lambda|  + \Lnorm{\infty}{f} + \Lnorm{\infty}{g})
 \Lnorm{\infty}{f- g}.
\end{equation*}
Employing this inequality we have the bound
\begin{eqnarray*}
\lefteqn{
\hspace{-3cm}
\Pt \Bigl(  \abs{R(\qt_{n},\lambda)-R(\qt, \lambda)} \geq \eta \lambda^2, \Lnorm{\infty}{\qt}\leq M | \lambda | , 
                                 \Lnorm{\infty}{\qt_{n}} \leq M  | \lambda |  \Bigr) } \\
     &  & \le 
 \Pt \Bigl( \Lnorm{\infty}{\qt-\qt_{n}}> \eta |\lambda|/(1+2M) \Bigr).
\end{eqnarray*}
Then using Borell's inequality once more, 
\begin{eqnarray}
\lefteqn{ \hspace{-4.3cm} \Pt  \Bigl( \abs{R(\qt_{n},\lambda)-R(\qt,\lambda)} 
\geq \eta \lambda^2, \Lnorm{\infty}{\qt}\leq M |\lambda| , 
\Lnorm{\infty}{\qt_{n}} \leq M |\lambda|  \Bigr)} \\
&  &  \leq C \exp{\left[- \frac{\eta^2 \lambda^2}{M^2 \sigma^2_{n}} \right] },
\nonumber
\label{ineq:No2}
\end{eqnarray}
for $|\lambda|$  large enough. Here
\[
\sigma^2_{n}= \sum_{\abs{k}>n} \lambda_k^2, 
\]
and we note for later that $ \sigma^2_n \downarrow 0$ as $n \uparrow \infty$.

\bigskip

\noindent
{\em Term 4}:
Denote by $I(\qt) = \Rprod{\qt}{{\bf K}^{-1} \qt}= \sum \lambda_k^{-2} \abs{c_k}^2$. 
From the study of the variational problem, Theorem \ref{thm:Lower} ,  
we have that, for any $\epsilon>0$,
\[
I(\qt) + R(\qt,\lambda) \geq \lambda^2 \times \frac{(1-\epsilon)}{K(0)}
\]
by choice of $\lambda \ll -1$. Thus, on the set $D(\gamma+2\eta)$ we may take 
\[
I(\qt) \geq \lambda^2 \frac{(1-\gamma-3\eta)}{K(0)} 
\]
while $\lambda \rightarrow - \infty$. It follows that
\begin{eqnarray}
\Pt \Bigl( \qt_{n} \in D(\gamma+2\eta)\Bigr) 
&\leq & \Pt \Bigl( I(\qt_{n}) \geq (1-\gamma-3\eta) \lambda^2 \} \Bigr)  \\
& \leq & \exp \Bigl[ -\frac{\lambda^2(1-\gamma-3\eta)(1-\eta)}{2 K(0)} \Bigr] 
\times  \tilde{E} \left[\exp \Bigl[{\frac{1-\eta}{2}I(\qt_{n})} \Bigr] \right] .
  \nonumber
\label{BoundRateFunction}
\end{eqnarray}
This  last factor is a simple Gaussian integral, in particular
\begin{eqnarray*}
\tilde{E} \left[\exp \Bigl[{\frac{1-\eta}{2}I(\qt_{n})} \Bigr] \right]  = 
\prod_{\abs{k}\leq n}\int_{\CC} \exp \left[{\frac{(1-\eta)\abs{g_k}^2}{2} }
  \right]   e^{-\frac{\abs{g_k}^2}{2} }\frac{d^2g_k}{{2\pi}}
= \eta^{-(2n+1)/2}.
\end{eqnarray*}
Together with (\ref{BoundRateFunction}) this produces
\begin{equation}
\Pt \Bigl( \qt_{n} \in D(\gamma+ 2 \eta \Bigr) \leq \eta^{-(2n+1)/2}
\exp \left[ -\frac{ \lambda^2(1-\gamma-3\eta)(1-\eta)}{2 K(0)} \right].
\label{ineq:No3}
\end{equation}
as the needed estimate on the fourth term.

\bigskip

Inequalities (\ref{ineq:No1}), (\ref{ineq:No2}) and (\ref{ineq:No3}) are now substituted 
into (\ref{ineq:DecompProbD}) to obtain
\begin{eqnarray*}
\Pt  \Bigl( H(\gamma,\eta) \Bigr) &\leq & 2C \exp{\left[- \frac{M^2 \lambda^2}{2 \sigma^2_{\infty}} \right] } + 
 C \exp{\left[- \frac{\eta^2 \lambda^2}{M^2 \sigma^2_{n}} \right] } \\
& & + \eta^{-(2n+1)/2}
\exp \left[ -\frac{ \lambda^2(1-\gamma-3\eta)(1-\eta)}{2 K(0)} \right].
\end{eqnarray*}
By choosing first $M$, then $n$, large  we may guarantee that 
the first two terms on the right hand side are negligible compared with the third.
Then 
taking logarithms, dividing by $\lambda^2$,  with $\lambda \rightarrow -\infty$ we
find that 
\begin{eqnarray*}
\limsup_{\lambda \rightarrow - \infty} 
\frac{1}{\lambda^2}\log\int_{H(\gamma,\eta)} e^{-\frac{R(\lambda,\qt)}{2}} d\Pt(\qt) 
& \le & -\frac{( 1-\gamma-3\eta)(1-\eta)}{2 K(0)}+ \frac{(-\gamma+\eta)}{2} \\
& \le &  - \frac{(1-\eta)}{2 K(0)}  + 2 \eta,
\end{eqnarray*}
for any $\gamma \in [0,1]$.  Finally recalling (\ref{ineq:Initial}) and (\ref{Cheby}), it also holds  
\[
\limsup_{\lambda \rightarrow  - \infty} \frac{1}{\lambda^2} \log f(\lambda)
\leq -\frac{( 1-\eta)}{2 K(0)}  + 2 \eta,
\]
and letting $\eta \downarrow 0$ completes the proof.

\bigskip
\noindent
{\bf Acknowledgements}   The research of B.R. was supported in part by NSF grant DMS-0505680. 
Also, most of this work was completed while J.R. visited the  Math. Dept. at CU Boulder; 
he is indebted to its members
for their hospitality.

\sc \bigskip \noindent 
Jos{\'e}  A. Ram{\'{\i}}rez \\
Department of Mathematics, \\ Universidad de Costa Rica, \\
San Jose  2060, Costa Rica. \\
{\tt  jaramirez@cariari.ucr.ac.cr}

\sc \bigskip \noindent Brian Rider \\ Department
of Mathematics,  \\ University of Colorado at Boulder, \\
UCB 395, Boulder, CO 80309. \\{\tt brider@euclid.colorado.edu}


\begin{thebibliography}{99}


\bibitem{AT} {\sc Adler, R. J. and Taylor, J. E.} (2006)
{\it Random Fields and Geometry}. To be published by Birkhauser.
({\em See iew3.technion.ac.il/\~{}radler/grf.pdf})

\bibitem{CM} {\sc Cambronero, S. and McKean, H. P.} (1999).
The ground state eigenvalue of Hill's equation with White Noise potential.
{\it Comm. Pure Appl. Math.}  {\bf 52}, 1277-1294.

\bibitem{CRR} {\sc Cambronero, S., Rider, B. and Ram\'{\i}rez, J.} (2006) 
{On the shape of the ground state eigenvalue density of a random Hill's equation.}
{\it Comm. Pure Appl. Math.}  {\bf 59}, 935-976.

\bibitem{DS}
{\sc Deuschel, J-D., and Stroock, D.W.} (1989)
{\it Large Deviations.} 
Academic Press, Boston.

\bibitem{FM}
{\sc Fukushima, M., and Nakao, S.} (1976/77)
On the spectra of the Schr\"odinger operator with
a white Gaussian noise potential.
{\it Z. Wahr. und Verw. Gabiete}
{\bf 37}  3, 267-274.

\bibitem{Hlp}
{\sc  Halperin, B.I.} (1965)
Green's functions for a particle in a one-dimensional random
potential.
{\it Phys. Rev. (2)} {\bf 139}, A104-A117.

\bibitem{LGP} 
{\sc Lifshits, I.M., Gredeskul, S. A. and Pastur, L.A.} (1988)
{\it Introduction to the theory of disordered systems}.
J. Wiley \& Sons, New York.

\bibitem{Mrk}
{\sc Merkl, F.} (2003)
Quenched asymptotics of the ground state energy
of random Schr\"odinger operators
with scaled Gibbsian potentials.
{\it Probab. Theory Relat. Fields}
{\bf 126}, 307-338.

\bibitem{M1}
{\sc McKean, H.P.} (1994)
A limit law for the groundstate of Hill's equation.
{\it J. Stat. Phys.} {\bf 74} no. 5-6, 1227-1232.


\bibitem{P1}
{\sc Pastur, L.A.} (1972)
The distribution of eigenvalues of the Schr\"odinger 
equation with a random potential.
{\it Functional Anal. Appl.} {\bf 6},  163-165.

\bibitem{PF}
{\sc Pastur, L.A.,  Figotin, A.} (1992)
{\it Spectra of random and almost periodic operators}.
Springer-Verlag, Belin \& Heidelberg.

 
\bibitem{RRV}
{\sc  Ram\'{\i}rez, J.,  Rider, B., and Vir\'ag, B.} (2006)
Beta ensembles, Stochastic Airy spectrum, and a diffusion.
{\it Preprint,  arXiv:math.PR:0607331.}

\bibitem{Sn}
{\sc Sznitman, A.-S.} (1998)
{\it Brownian motion, obstacles and random media.}
Springer Monographs in Mathematics, Berlin-Heidelberg.



\end{thebibliography}
\end{document}